\documentclass[A4,12pt]{article}

\usepackage{lineno,hyperref}
\usepackage{subfigure}
\usepackage{graphicx}
\usepackage{mathrsfs}
\usepackage{algorithm}
\usepackage{algpseudocode}
\algtext*{EndProcedure}
\begin{document}

\title{A numerical model based on the curvilinear coordinate system for the MAC method simplified \\} 

\author{\\ Eliandro Rodrigues Cirilo, Alessandra Negrini Dalla Barba,\\ Neyva Maria Lopes Romeiro, Paulo Laerte Natti \\ \\
Universidade Estadual de Londrina,\\ Rod. Celso Garcia Cid, PR-445, km 380, \\ 86051-990, Londrina - PR, Brazil}

\maketitle

\begin{abstract}
In this paper we developed a numerical methodology to study some incompressible fluid flows without free surface, using the curvilinear coordinate system and whose edge geometry is constructed via parametrized spline. First, we discussed the representation of the Navier-Stokes and continuity equations on the curvilinear coordinate system, along with the auxiliary conditions. 
Then, we presented the numerical method -- a simplified version of MAC (\textit{Marker and Cell}) method -- along with the discretization of the governing equations, which is carried out using the finite differences method and the implementation of the FOU (\textit{First Order Upwind}) scheme. 
Finally, we applied the numerical methodology to the parallel plates problem, lid-driven cavity problem and atherosclerosis problem, and then we compare the results obtained with those presented in the literature.
\end{abstract}

\noindent \textbf{keyword}: finite difference method, simplified MAC method, curvilinear coordinate system, parallel plates problem, lid-driven cavity problem, atherosclerosis problem 

\section{Introduction}

In order to find numerical solutions to a differential equations it is indispensable to approximate the 
equations by a mathematical method that can be programmed in the computers.
There are some techniques used for this purpose for example - finite differences \cite{griebel98},
finite element \cite{johnson_1987} , finite volume \cite{leveque_2002}, SPH
(Smoothed Particles Hydrodynamics) \cite{Monaghan_2005}, among others. We use finite differences in this work. 

In general, the geometry which represent the physical domain of the problems are irregular forms. 
That is why we chose to study the 
curvilinear coordinate system. This system allows to generate mesh that coincides with the boundary of 
the domain in the applied problems \cite{Pardo_2012,Romeiro_2011,Romeiro17}. 
In addition, we can to represent the differential equations, 
to be solved, in this system \cite{thompson85,maliska13}.

When the partial differential equations are considered in the formulation of the mathematical model
of a specific problem it is important to know the terms that compose these equations. The terms can be temporals, convectives, 
diffusives and others. 
The proper way to 
discretize such terms, respecting inherent physical of the problem, is necessary for the success of obtaining the numerical solution. 
In particular, the development of approximate methods of convective terms have been studied by many research in recent years \cite{ferreira12}. The accuracy 
of the results obtained from the numerical solutions is directly influenced by the choice of convection scheme.
The \textit{upwind} scheme can be used in this purpose. This approach is done according with the sign of local 
convection speed. The \textit{upwind} schemes are classified as first order or high order. In the first category we can 
mention the FOU scheme (\textit{First Order Upwinding}) \cite{courant52}, which is stable unconditionally and produces a diffusive 
character which usually
smooths the solution. Among the high-order we can mention the SOU (\textit{Second Order Upwinding}) \cite{Price1966},
QUICK (\textit{Quadratic Upstream Interpolation for Convective Kinematics}) \cite{leonard_1979}, which contributing to 
increase the accuracy of the numerical method, but introduce non-physical 
oscillations that can compromise the convergence.
Finally, the third-order accurate upwind compact finite difference schemes,
which allowed to obtain accurate numerical results for the benchmark flow problems \cite{Shah_2012}.

In this work the objective is propose a versatile numerical method to study some incompressible fluids flow without
surface free, in the curvilinear coordinate system, using the \textit{upwind} scheme FOU to approximate the convective terms,
with the board interpolated by parametrized spline, which allows better simulate a greater amount of complex problems.
This numerical method is applied to the study of three problems of incompressible fluid flows: parallel plates problem, lid-driven cavity problem, and atherosclerosis problem. In addition, the work shows that our numerical methodology accurately reproduces the results observed in the literature of these three problems studied.

\section{Curvilinear governing equations}

In many cases, when we want to get analytical solutions to the problems under study, many simplifications are necessary, and this
can depart us of the reality. When our goal is to study a certain physical phenomenon more realistically, we need initially to
model the physics problem.
As usually, the equations obtained in this modeling process has no analytical solution, numerical methods are then used to obtain the
solution of the problem studied.
Since it isn't possible to obtain numerical solutions on the
continuous region, as this involves infinite points, we first discretize the domain, that is, divide it in points and only in these is that 
we'll find the solution of the problem. This set of points is called mesh. 
It is important that they are properly distributed so that 
the numerical solution represents satisfactorily the gradients of interest in the problem studied.

Generally, problems of everyday life are not evaluated in rectangular coordinates, because are represented by irregular geometries. 
Then, the curvilinear coordinate system it is more appropriate. Its main function is the representation of complex 
geometries, in these cases the cartesian coordinate system leads to a poor fit of the border, since the
physical domain doesn't match the domain mesh \cite{thompson85}.

By means of a transformation (which may be numeric) between the cartesian coordinate system $(x, y)$ and the curvilinear coordinate 
system $(\xi, \eta)$, it is possible to map the domain written in the $(x,y)$ system to other written on a regular geometry $(\xi, \eta)$.
The $(x, y)$ system is termed physical domain, while the $(\xi, \eta)$ is called computational domain. 

\begin{figure}[!ht]
\begin{center}
\caption{Representation of the $R$ region}
\subfigure[Physical domain]{\includegraphics[height=4.4cm,width=5.4cm]{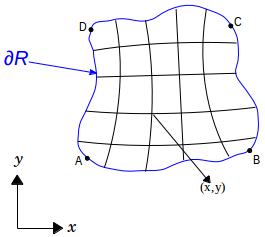}}
   \label{domfis}\qquad \qquad
\subfigure[Computational domain]{\includegraphics[height=4.4cm,width=5.4cm]{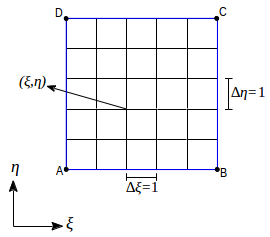}
    \label{domcomput}}
%\caption{Representation of the $R$ region}
\end{center}
\center{{\bf{Source:}} The Authors}
\label{regiaor}
\end{figure}

For example, in the two-dimensional case, the computational domain is taken in rectangular form, regardless of the
physical geometry, Figs. 1(a) and 1(b). 
Thus, the physical domain is transformed so the computational domain is always represented by a rectangle.
For convenience it is assumed elementary volumes with unitary dimensions, i.e.,
$\Delta \xi = \Delta \eta = 1$.
So, as much as the coordinate lines take arbitrary spacing on the physical plane, in the computational
dimensions are fixed.

The coordinates of an arbitrary point in the curvilinear coordinate system $(\xi,\eta)$ are related with the cartesian coordinate system 
$(x, y)$ for transformation equations as follows
\begin{eqnarray*}
\label{xietatau}
 \xi = \xi(x,y,t); \quad 
 \eta = \eta(x,y,t); \quad
 \tau = \tau(t)
\end{eqnarray*}
\noindent where $\tau = t$, because we are not admitting mesh movement. Therefore, the metric transformation are given by 
 \begin{eqnarray}
 \label{mettrans}
  \frac{\partial{\xi}}{\partial{x}} = J \frac{\partial{y}}{\partial{\eta}}; \quad
  \frac{\partial{\xi}}{\partial{y}} = - J \frac{\partial{x}}{\partial{\eta}}; \quad
  \frac{\partial{\eta}}{\partial{x}} = - J \frac{\partial{y}}{\partial{\xi}}; \quad
  \frac{\partial{\eta}}{\partial{y}} = J \frac{\partial{x}}{\partial{\xi}} 
 \end{eqnarray}
\noindent where
\begin{eqnarray}
\label{jacob}
  J = \left( \frac{\partial{x}}{\partial{\xi}} \frac{\partial{y}}{\partial{\eta}} 
              - \frac{\partial{x}}{\partial{\eta}} \frac{\partial{y}}{\partial{\xi}} \right)^{-1}
\end{eqnarray}
\noindent is the jacobian of the transformation \cite{thompson85,maliska13}.

The metrics enabled mapping of the physical domain for the computational, implying in the perform
of the required geometric compensations.
In this work we build the mesh as detailed in \cite{cirilo06,saita17} and the board was obtained via the 
parametrized Spline method.
Moreover, mathematical compensations is also carried out, by chain rule,
in the governing equations of the studied physical problem. Therefore, the conformity between the computational
mesh and the governing equations allows adequately handle the computational numerical simulation.

\begin{figure}[!ht]
\centering
\caption{Cells}
\subfigure[Structure of the mesh and its nomenclature]{
\includegraphics[height=5cm]{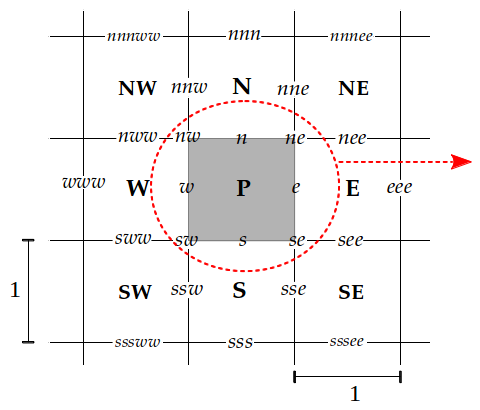}
\label{malhansew2}
}
\subfigure[Pressure and velocity storage]{
\includegraphics[height=4cm]{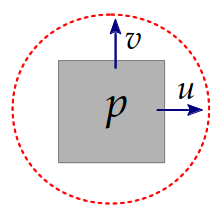}
\label{armazenamento}
}
%\caption{Cells}
\center{{\bf{Source:}} The Authors}
\label{celulas}
\end{figure} 

So we can solve the problems of interest. First, it is need to label the points in the mesh such a way that calculations are performed 
correctly by the numerical method.
Considering Fig. \ref{celulas}, we have that the mesh is composed by cells where the topological relationships are arranged and 
labeled as the cardinal points. The labels $P$, $E$, $W$, $N$, $S$, $NE$, 
$SE$, $NW$, $SW$ mean center, east, west, north, south, northeast, southeast, northwest, southwest, respectively. The abbreviations in 
tiny, 
positioned on the sides, are cardinal changes from the center of the cell labeled by $P$. 

Assuming by hypothesis a laminar flow, newtonian, isothermal and incompressible, 
in two-dimensional, without the existence of the source terms, 
then the Navier-Stokes equations described in the curvilinear coordinate system are written by
\begin{eqnarray}
\label{coordgeneru}
       \underbrace{ \frac{\partial{}}{\partial{\tau}} \left( \frac{\rho u}{J} \right) }_{\mbox{ temporal term }}
        & + & \underbrace{ \frac{\partial{}}{\partial{\xi}} \left( \rho U u \right) 
        + \frac{\partial{}}{\partial{\eta}} \left( \rho V u \right) }_{\mbox{convective term}}
       = \underbrace{ \left[ \frac{\partial{p}}{\partial{\eta}} \frac{\partial{y}}{\partial{\xi}}
       - \frac{\partial{p}}{\partial{\xi}} \frac{\partial{y}}{\partial{\eta}} \right] }_{\mbox{pressure term}} \\ 
       & + & \underbrace{ \mu \left[ \frac{\partial{}}{\partial{\xi}} \left( J  
       \left( \alpha \frac{\partial{u}}{\partial{\xi}} - \beta \frac{\partial{u}}{\partial{\eta}} \right) \right)
       + \frac{\partial{}}{\partial{\eta}} \left( J \left( \gamma \frac{\partial{u}}{\partial{\eta}} 
       - \beta \frac{\partial{u}}{\partial{\xi}} \right) \right) \right] }_{\mbox{diffusive term}} \nonumber
\end{eqnarray}

\begin{eqnarray}
\label{coordgenerv}
       \underbrace{ \frac{\partial{}}{\partial{\tau}} \left( \frac{\rho v}{J} \right) }_{\mbox{temporal term}}
       & + & \underbrace{ \frac{\partial{}}{\partial{\xi}} \left( \rho U v \right) 
       + \frac{\partial{}}{\partial{\eta}} \left( \rho V v \right) }_{\mbox{convective term}}
       = \underbrace{ \left[ \frac{\partial{p}}{\partial{\xi}} \frac{\partial{x}}{\partial{\eta}}
       - \frac{\partial{p}}{\partial{\eta}} \frac{\partial{x}}{\partial{\xi}} \right] }_{\mbox{pressure term}} \\ 
       & + & \underbrace{ \mu \left[ \frac{\partial{}}{\partial{\xi}} \left( J
       \left( \alpha \frac{\partial{v}}{\partial{\xi}} - \beta \frac{\partial{v}}{\partial{\eta}} \right) \right)
       + \frac{\partial{}}{\partial{\eta}} \left( J \left( \gamma \frac{\partial{v}}{\partial{\eta}} 
       - \beta \frac{\partial{v}}{\partial{\xi}} \right) \right) \right] }_{\mbox{diffusive term}} \nonumber
\end{eqnarray}

\noindent where the viscosity $\mu$ is constant. For details of how to describe this equations in these coordinate system see 
\cite{maliska_1983,thompson76}.
The terms $U$ and $V$ in (\ref{coordgeneru}) and (\ref{coordgenerv}) are named contravariant components of the velocity vector, 
these are normal to $\xi$ and $\eta$ lines respectively, and defined as 
$
\displaystyle
\label{contrav}
  U = \frac{1}{J} \left( u \frac{\partial{\xi}}{\partial{x}} 
       + v \frac{\partial{\xi}}{\partial{y}} \right)
$
and
$
\displaystyle
  V = \frac{1}{J} \left( u \frac{\partial{\eta}}{\partial{x}} 
       + v \frac{\partial{\eta}}{\partial{y}} \right) \nonumber
$. As the incompressibility hypothesis was assumed, the continuity equation is such that
\begin{eqnarray}
\label{conservmassagen}
   \frac{\partial{U}}{\partial{\xi}} + \frac{\partial{V}}{\partial{\eta}} = 0
\end{eqnarray}
\noindent therefore the governing equations of our computational model, in curvilinear coordinates,
are given by (\ref{coordgeneru}), (\ref{coordgenerv}) and (\ref{conservmassagen}). 

The choice of auxiliary conditions is critical to the formulation of any problem described by differential equations 
\cite{griebel98,fortuna12}. As 
the Navier-Stokes equations may be used to describe the flow in many situations, it is important to properly set the initial and boundary
conditions of the problem under study.
For the initial conditions is taken a velocity field that satisfies continuity equation. The boundary conditions used in the resolution 
of problems considered can be classified in no-slip and impermeability condition (CNEI), free-slip condition (CLES), 
prescribed injection condition (CIPR) and continuous ejection condition (CECO).
The details about the numerical considerations made on the governing equations and auxiliary conditions are detailed as follows.

\section{Numerical modeling}

For the numerical solution of the fluid equations has been chosen
the MAC (\textit{Marker and cell}) method formulated by Harlow and Welch \cite{Amsden_1970}. This method has the 
advantage of permitting the simulation of different types of flow, in the cartesian coordinate system, with and 
without free surface \cite{McKee_2008,tome_2014} and multiphase flows \cite{santos_2012}.

In \cite{Patil_2009} the authors only describe the procedure of how to solve the mechanical equations of fluid via the
curvilinear coordinates system. However, what innovates in this study is that we describe the numerical procedure
considering: (1) the governing equations in dimensional form, (2) to an irregular geometry 
whose board is obtained by parametric Spline, and (3) the system of linear equations for the pressure resolution is solved 
by Gauss-Seidel methodology. In this work has been considered a MAC methodology to confined flow.
With this action has been disregarded any moving marker particles deduction associated with the method, which results 
in simplification of the numerical calculations.

Furthermore it is considered displaced mesh, where the unknowns are stored in different positions. 
The displaced storage for the components of the velocity 
vector and pressure has a positive impact on numerical calculation, due to the fact that reduce numerical instability 
\cite{fortuna12,Amsden_1970}.
As shown in Figs. \ref{malhansew2} and \ref{armazenamento}, the pressure ($p$) is located in the center of the cell and the components of velocity vector 
($u$ and $v$) in the centers of the faces. Note that $p$ refers to the pressure and $P$ the center of the cell.
Rewriting the equations (\ref{coordgeneru}) and (\ref{coordgenerv}) as follows
\begin{eqnarray*}
\label{equgen1}
\qquad \frac{\partial{}}{\partial{\tau}} \left( \frac{u}{J} \right) =
     - { \mathscr{C}(u) } + \frac{1}{\rho} { \mathscr{P}^{u} } + \nu { \mathscr{V}(u) }
\quad ; \quad
\label{eqvgen1}
\frac{\partial{}}{\partial{\tau}} \left( \frac{v}{J} \right) =
     - { \mathscr{C}(v) } + \frac{1}{\rho} { \mathscr{P}^{v} } + \nu { \mathscr{V}(v) }
\end{eqnarray*}

\noindent where $\displaystyle\frac{\partial{}}{\partial{\tau}} \left( \frac{u}{J} \right)$ and 
$\displaystyle\frac{\partial{}}{\partial{\tau}} \left( \frac{v}{J} \right)$ are the temporal terms, $\mathscr{C}(u)$ and $\mathscr{C}(v)$
the convective terms, $\mathscr{P}^u$ and $\mathscr{P}^v$ the pressure terms and $\mathscr{V}(u)$ and $\mathscr{V}(v)$ the diffusive
terms and $\nu = \mu/\rho$ the viscosity term.

For the temporal term were
done first-order approximation. In the case of convective terms apply approaches \textit{upwind} type. The pressure terms are 
approximated by central differences. Finally, the diffusive terms are approximate by central differences, too.
The discretization of the equations on the face $e$, for every cell within the domain, in the time level
$k$, is written as 
\begin{eqnarray}
\label{discretiz1}
\displaystyle
\left. \frac{\partial{}}{\partial{\tau}} \left( \frac{u}{J} \right) \right|_{e}^{k}
      = \frac{1}{J_{e}} \left( \frac{\left. u \right|_{e}^{k+1} - \left. u \right|_{e}^{k}}{\Delta \tau} \right)
      = - \mathscr{C}(u)|_{e}^{k} + \frac{1}{\rho} \mathscr{P}^{u}|_{e}^{k+1} + \nu \mathscr{V}(u)|_{e}^{k} \nonumber
\end{eqnarray}
\noindent analogously
\begin{eqnarray}
\label{discretiz2}
\displaystyle
\left. \frac{\partial{}}{\partial{\tau}} \left( \frac{v}{J} \right) \right|_{n}^{k}
      = \frac{1}{J_{n}} \left( \frac{\left. v \right|_{n}^{k+1} - \left. v \right|_{n}^{k}}{\Delta \tau} \right)
      = - \mathscr{C}(v)|_{n}^{k} + \frac{1}{\rho} \mathscr{P}^{v}|_{n}^{k+1} + \nu \mathscr{V}(v)|_{n}^{k} \nonumber
\end{eqnarray}

Setting
$\displaystyle
\label{eqFek}
F|_{e}^{k} = J_{e} \Delta \tau \left[ - \mathscr{C}(u)|_{e}^{k} + \nu \mathscr{V}(u)|_{e}^{k} \right] + u|_{e}^{k}
$ the cartesian component $u$ of the velocity vector is expressed as
\begin{eqnarray}
\label{equek1}
  u|_{e}^{k+1} = F|_{e}^{k} + \frac{J_{e} \Delta \tau}{\rho} \mathscr{P}^{u}|_{e}^{k+1}
\end{eqnarray}

\noindent and denoting
$\displaystyle
\label{eqGnk}
  G|_{n}^{k} = J_{n} \Delta \tau \left[ - \mathscr{C}(v)|_{n}^{k} + \nu \mathscr{V}(v)|_{n}^{k} \right] + v|_{n}^{k}
$ we have
\begin{eqnarray}
\label{eqvnk1}
  v|_{n}^{k+1} = G|_{n}^{k} + \frac{J_{n} \Delta \tau}{\rho} \mathscr{P}^{v}|_{n}^{k+1}
\end{eqnarray}

One of the great difficulties in numerically solving the Navier-Stokes is about the discretization of the convective terms, because these terms are not linear. We use the FOU scheme for these terms.
Of course this is not the most appropriate scheme to very low/high Reynolds number, but in this study we want to show that the numerical model is appropriate.

The discretization of the $\mathscr{C}(u)$ term on the face $e$ of the cell centered in the cardinal point $P$, in a time level $k$,
using central differences, can be given as
\begin{eqnarray*}
\displaystyle \nonumber \left. \mathscr{C}(u) \right|_{e}^{k} 
     \approx \overline{U}_{E}^{k} u_{E}^{k} - \overline{U}_{P}^{k} u_{P}^{k}
        + \overline{V}_{ne}^{k} u_{ne}^{k} - \overline{V}_{se}^{k} u_{se}^{k}
\end{eqnarray*}

\noindent with convection velocities calculated by arithmetic mean, and applying FOU scheme to approximate $u$ component 
in the corresponding face.
On the other hand, the discretization of the $\mathscr{C}(v)$ on the face $n$, in time level $k$, is given by
\begin{eqnarray*}
\displaystyle
 \nonumber \left. \mathscr{C}(v) \right|_{n}^{k} \approx  \overline{U}_{ne}^{k} v_{ne}^{k} - \overline{U}_{nw}^{k} v_{nw}^{k}
                     + \overline{V}_{N}^{k} v_{N}^{k} - \overline{V}_{P}^{k} v_{P}^{k}
\end{eqnarray*}
\noindent with the convection velocities calculated by arithmetic mean, and approaches with the \textit{upwind} scheme, too. 
For more details see \cite{fortuna12}. 
The discretization of the pressure terms on the faces $e$ and $n$, in time level $k+1$, are as follows
\begin{eqnarray*}
\displaystyle
 \nonumber \left. \left( \frac{1}{\rho} \mathscr{P}^u \right) \right|_{e}^{k+1} 
      \approx \frac{1}{\rho} \left[ \left( p_{ne}^{k+1} - p_{se}^{k+1} \right) \left. \frac{\partial{y}}{\partial{\xi}} \right|_{e} 
     - \left( p_{E}^{k+1} - p_{P}^{k+1} \right) \left. \frac{\partial{y}}{\partial{\eta}} \right|_{e} \right]
\end{eqnarray*}
\noindent and
\begin{eqnarray*}
 \nonumber \left. \left( \frac{1}{\rho} \mathscr{P}^v \right) \right|_{n}^{k+1} 
 \approx \frac{1}{\rho} \left[ \left( p_{ne}^{k+1} - p_{nw}^{k+1} \right) \left. \frac{\partial{x}}{\partial{\eta}} \right|_{n} 
     - \left( p_{N}^{k+1} - p_{P}^{k+1} \right) \left. \frac{\partial{x}}{\partial{\xi}} \right|_{n} \right]
\end{eqnarray*}

Finally, the discretization of the diffusive terms in faces $e$ and $n$, in time level $k$, can be given as
\begin{eqnarray*}
\displaystyle
\label{exp5}
 \nonumber \left. \left( \nu \mathscr{V}(u) \right) \right|_{e}^{k}
    = \nu_{e} \left. \left[ \frac{\partial{}}{\partial{\xi}} \left( J \left( \alpha \frac{\partial{u}}{\partial{\xi}}
         - \beta \frac{\partial{u}}{\partial{\eta}} \right) \right) \right] \right|_{e}^{k} 
         + \nu_{e} \left. \left[ \frac{\partial{}}{\partial{\eta}} \left( J \left( \gamma \frac{\partial{u}}{\partial{\eta}}
         - \beta \frac{\partial{u}}{\partial{\xi}} \right) \right) \right]\right|_{e}^{k}
\end{eqnarray*}
\begin{eqnarray*}
\label{exp6}
 \nonumber \left. \left( \nu \mathscr{V}(v) \right)\right|_{n}^{k}
     = \nu_{n} \left. \left[ \frac{\partial{}}{\partial{\xi}} \left( J \left( \alpha \frac{\partial{v}}{\partial{\xi}}
         - \beta \frac{\partial{v}}{\partial{\eta}} \right) \right) \right]\right|_{n}^{k} 
         + \nu_{n} \left. \left[ \frac{\partial{}}{\partial{\eta}} \left( J \left( \gamma \frac{\partial{v}}{\partial{\eta}}
         - \beta \frac{\partial{v}}{\partial{\xi}} \right) \right) \right]\right|_{n}^{k}
\end{eqnarray*}
\noindent with
\begin{eqnarray*}
\label{exp1}
\displaystyle
\nonumber \left. \frac{\partial{}}{\partial{\xi}} \left(J \left( \alpha \frac{\partial{u}}{\partial{\xi}} 
            - \beta \frac{\partial{u}}{\partial{\eta}} \right) \right) \right|_{e}^{k} 
   & \approx & \left. (J \alpha)\right|_{E} \left( u_{eee}^{k} - u_{e}^{k} \right) - \left. (J \beta)\right|_{E} \left( u_{nee}^{k} - u_{see}^{k} \right) \\ 
            - \left. (J \alpha)\right|_{P} \left( u_{e}^{k} - u_{w}^{k} \right) 
            &+& \left. (J \beta)\right|_{P} \left( u_{n}^{k} - u_{s}^{k} \right)
\\ \nonumber \\             
\label{exp2}
\nonumber \frac{\partial{}}{\partial{\eta}} \left(J \left( \gamma \frac{\partial{u}}{\partial{\eta}} 
            - \beta \frac{\partial{u}}{\partial{\xi}} \right) \right)|_{e}^{k}
            & \approx & (J \gamma)|_{ne}^{k} \left( u_{nne}^{k} - u_{e}^{k} \right)
            - (J \beta)|_{ne}^{k} \left( u_{nee}^{k} - u_{n}^{k} \right) \\
             -  (J \gamma)|_{se}^{k} \left( u_{e}^{k} - u_{sse}^{k} \right)
            &+& (J \beta)|_{se}^{k} \left( u_{see}^{k} - u_{s}^{k} \right) 
\\ \nonumber \\
\label{exp3}
\nonumber \frac{\partial{}}{\partial{\xi}} \left(J \left( \alpha \frac{\partial{v}}{\partial{\xi}} 
            - \beta \frac{\partial{v}}{\partial{\eta}} \right) \right)|_{n}^{k}
            & \approx & (J \alpha)|_{ne}^{k} \left( v_{nee}^{k} - v_{n}^{k} \right)
            - (J \beta)|_{ne}^{k} \left( v_{nne}^{k} - v_{e}^{k} \right) \\ 
            - (J \alpha)|_{nw}^{k} \left( v_{n}^{k} - v_{nww}^{k} \right)
            &+& (J \beta)|_{nw}^{k} \left( v_{nnw}^{k} - v_{w}^{k} \right) 
\\ \nonumber \\
\label{exp4}
\nonumber \frac{\partial{}}{\partial{\eta}} \left(J \left( \gamma \frac{\partial{v}}{\partial{\eta}} 
            - \beta \frac{\partial{v}}{\partial{\xi}} \right) \right)|_{n}^{k}
            & \approx & (J \gamma)|_{N}^{k} \left( v_{nnn}^{k} - v_{n}^{k} \right)
            - (J \beta)|_{N}^{k} \left( v_{nne}^{k} - v_{nnw}^{k} \right) \\ 
            - (J \gamma)|_{P}^{k} \left( v_{n}^{k} - v_{s}^{k} \right)
            &+& (J \beta)|_{P}^{k} \left( v_{e}^{k} - v_{w}^{k} \right)
\end{eqnarray*}

Since the contravariant components $U$ and $V$ can be rewritten as 
$\displaystyle U = u \frac{\partial{y}}{\partial{\eta}} - v \frac{\partial{x}}{\partial{\eta}}$
and
$\displaystyle V = - u \frac{\partial{y}}{\partial{\xi}} + v \frac{\partial{x}}{\partial{\xi}}$,
from the expressions (\ref{equek1}) and (\ref{eqvnk1}) -- and similar for the other faces -- we get the
compact form for contravariant components
\begin{eqnarray}
\label{uek1}
U|_{e}^{k+1} &\!\!\!\!=\!\!\!\! & F|_{e}^{k} \left. \frac{\partial{y}}{\partial{\eta}}\right|_{e} \!\!\!\! - G|_{e}^{k} \left. \frac{\partial{x}}{\partial{\eta}}\right|_{e}
      + \frac{J_{e} \Delta \tau}{\rho} \left\{ - \left. \frac{\partial{p}}{\partial{\xi}}\right|_{e}^{k+1} \!\!\!\! \alpha|_{e} 
      + \left. \frac{\partial{p}}{\partial{\eta}}\right|_{e}^{k+1} \!\!\!\! \beta|_{e} \right\} \\ \nonumber \\
\label{uwk1}
U|_{w}^{k+1} &\!\!\!\!=\!\!\!\! & F|_{w}^{k} \left. \frac{\partial{y}}{\partial{\eta}}\right|_{w} \!\!\!\! - G|_{w}^{k} \left. \frac{\partial{x}}{\partial{\eta}}\right|_{w}
      + \frac{J_{w} \Delta \tau}{\rho} \left\{ - \left. \frac{\partial{p}}{\partial{\xi}}\right|_{w}^{k+1} \!\!\!\! \alpha|_{w} 
      + \left. \frac{\partial{p}}{\partial{\eta}}\right|_{w}^{k+1} \!\!\!\! \beta|_{w} \right\} \\ \nonumber \\
\label{vnk1}
V|_{n}^{k+1} &\!\!\!\!=\!\!\!\! & - F|_{n}^{k} \left. \frac{\partial{y}}{\partial{\xi}}\right|_{n} \!\!\!\! + G|_{n}^{k} \left. \frac{\partial{x}}{\partial{\xi}}\right|_{n}
      + \frac{J_{n} \Delta \tau}{\rho} \left\{ \left. \frac{\partial{p}}{\partial{\xi}}\right|_{n}^{k+1} \!\!\!\! \beta|_{n} 
      - \left. \frac{\partial{p}}{\partial{\eta}}\right|_{n}^{k+1} \!\!\!\! \gamma|_{n} \right\} \\ \nonumber \\
\label{vsk1}
V|_{s}^{k+1} &\!\!\!\!=\!\!\!\! & - F|_{s}^{k} \left. \frac{\partial{y}}{\partial{\xi}}\right|_{s} \!\!\!\! + G|_{s}^{k} \left. \frac{\partial{x}}{\partial{\xi}}\right|_{s}
      + \frac{J_{s} \Delta \tau}{\rho} \left\{ \left. \frac{\partial{p}}{\partial{\xi}}\right|_{s}^{k+1} \!\!\!\! \beta|_{s} 
      - \left. \frac{\partial{p}}{\partial{\eta}}\right|_{s}^{k+1} \!\!\!\! \gamma|_{s} \right\}
\end{eqnarray}

Approaching by central difference scheme in P cardinal point, in time level k + 1, the continuity equation gives us 
the following expression 
\begin{eqnarray}
\label{ueuwvnvs}
\displaystyle
\left. \frac{\partial{U}}{\partial{\xi}}\right|_{P}^{k+1} + \left. \frac{\partial{V}}{\partial{\eta}}\right|_{P}^{k+1} = 0
\quad \Rightarrow \quad  U|_{e}^{k+1} - U|_{w}^{k+1} + V|_{n}^{k+1} - V|_{s}^{k+1} = 0
\end{eqnarray}

\noindent Replacing the equations (\ref{uek1}) to (\ref{vsk1}) in (\ref{ueuwvnvs}) and grouping like terms in the same side of equality, we find
\begin{eqnarray}
\label{pressnumerico}
% \displaystyle
 \nonumber J_{e} \left\{ - \left. \frac{\partial{p}}{\partial{\xi}}\right|_{e}^{k+1} \alpha|_{e} 
           + \left. \frac{\partial{p}}{\partial{\eta}}\right|_{e}^{k+1} \beta|_{e} \right\} 
 + J_{w} \left\{ \left. \frac{\partial{p}}{\partial{\xi}}\right|_{w}^{k+1} \alpha|_{w} 
           - \left. \frac{\partial{p}}{\partial{\eta}}\right|_{w}^{k+1} \beta|_{w} \right\}  \\ \nonumber \\
 + J_{n} \left\{ \left. \frac{\partial{p}}{\partial{\xi}}\right|_{n}^{k+1} \beta|_{n} 
          - \left. \frac{\partial{p}}{\partial{\eta}}\right|_{n}^{k+1} \gamma|_{n} \right\} \nonumber
 + J_{s} \left\{ - \left. \frac{\partial{p}}{\partial{\xi}}\right|_{s}^{k+1} \beta|_{s} 
          + \left. \frac{\partial{p}}{\partial{\eta}}\right|_{s}^{k+1} \gamma|_{s} \right\} \\ \nonumber \\
= \frac{\rho}{\Delta \tau} \left\{ 
- F|_{e}^{k} \left.\frac{\partial{y}}{\partial{\eta}}\right|_{e} + G|_{e}^{k} \left.\frac{\partial{x}}{\partial{\eta}}\right|_{e}
 + F|_{w}^{k} \left.\frac{\partial{y}}{\partial{\eta}}\right|_{w} - G|_{w}^{k} \left.\frac{\partial{x}}{\partial{\eta}}\right|_{w}\right.  \nonumber \\
 \left.+ F|_{n}^{k} \left.\frac{\partial{y}}{\partial{\xi}}\right|_{n} - G|_{n}^{k} \left.\frac{\partial{x}}{\partial{\xi}}\right|_{n} 
 - F|_{s}^{k} \left.\frac{\partial{y}}{\partial{\xi}}\right|_{s} + G|_{s}^{k} \left.\frac{\partial{x}}{\partial{\xi}}\right|_{s} \right\}
\end{eqnarray}
\noindent The equation (\ref{pressnumerico}) is the pressure evolution equation, that satisfies the continuity equation 
(\ref{conservmassagen}).

\begin{figure}[!ht]
\centering
\caption{Possible configurations to the boundary conditions}
\subfigure[Case 1]{\includegraphics[height=1.4cm]{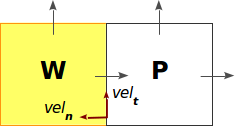}
\label{cnei1}}
\quad 
\subfigure[Case 2]{\includegraphics[height=1.4cm]{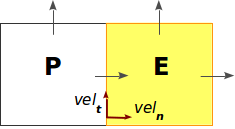}
\label{cnei2}}
\quad 
\subfigure[Case 3]{\includegraphics[height=2.6cm]{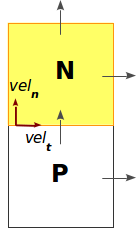}
\label{cnei3}}
\quad 
\subfigure[Case 4]{\includegraphics[height=2.6cm]{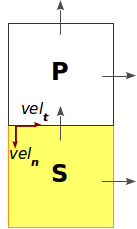}
\label{cnei4}}
%\caption{Possible configurations to the boundary conditions}
\label{condcnei}
\center{{\bf{Source:}} The Authors}
\end{figure}

The initial conditions considered should satisfy the continuity equation. We performed simulations always starting from a velocity 
field and pressure in a state of quiescence. 
We understand the quiescence state as the one with null velocity and pressure fields.
The pressure boundary condition is taken as $\displaystyle\frac{\partial{p}}{\partial{n}} = 0$.
On the boundary conditions, whatever the two-dimensional problem under study, there are four kinds of configurations between the 
contour and interior cells in the computational domain.
Denoting $vel_t$, $vel_n$ as the tangential and normal velocities at the border of the cells, $vel_{\bullet}$ (prescribed 
velocity), $vel_I$ (prescribed injection velocity), $vel_E$ (prescribed ejection velocity), 
then the four possible configurations are those shown in Fig. \ref {condcnei}.
\vspace{-0.5cm}

\vspace{-0.2cm}

\begin{description}
\item[Case 1:] The non-slip and impermeability condition (CNEI) is defined by the following expressions 
$ vel_t = 0 \mbox{ and } vel_n = 0 $, so
  \begin{eqnarray*}
      \mbox{\textbf{(a)}} \qquad
      vel_t = \frac{v|_{nww} + v|_{n}}{2} = 0   \qquad &\Rightarrow& \qquad    v|_{nww} = - v|_{n} \\ 
      vel_n = u|_{w} = 0   \qquad &\Rightarrow& \qquad u|_{w} = 0 \\
      \mbox{\textbf{(b)}} \qquad
      vel_t = \frac{v|_{n} + v|_{nee}}{2} = 0   \qquad &\Rightarrow& \qquad    v|_{nee} = - v|_{n} \\ 
      vel_n = u|_{e} = 0   \qquad &\Rightarrow& \qquad u|_{e} = 0 \\
      \mbox{\textbf{(c)}} \qquad
      vel_t = \frac{u|_{e} + u|_{nne}}{2} = 0   \qquad &\Rightarrow& \qquad    u|_{nne} = - u|_{e} \\ 
      vel_n = v|_{n} = 0   \qquad &\Rightarrow& \qquad v|_{n} = 0 \\
      \mbox{\textbf{(d)}} \qquad
      vel_t = \frac{u|_{e} + u|_{sse}}{2} = 0   \qquad &\Rightarrow& \qquad    u|_{sse} = - u|_{e} \\ 
      vel_n = v|_{s} = 0   \qquad &\Rightarrow& \qquad v|_{s} = 0
   \end{eqnarray*}

\item[Case 2:] The free-slip condition (CLES) is defined by the following expressions
$ vel_t = vel_{\bullet} \mbox{ and } vel_n = 0$, 
    \begin{eqnarray*}
      \mbox{\textbf{(a)}} \qquad
      vel_t = \frac{v|_{nww} + v|_{n}}{2} = vel_{\bullet}   \qquad &\Rightarrow& \qquad  v|_{nww} = 2 vel_{\bullet} - v|_{n} \\ 
      vel_n = u|_{w} = 0   \qquad &\Rightarrow& \qquad u|_{w} = 0 \\
      \mbox{\textbf{(b)}} \qquad
      vel_t = \frac{v|_{n} + v|_{nee}}{2} = vel_{\bullet}   \qquad &\Rightarrow& \qquad    v|_{nee} = 2 vel_{\bullet} - v|_{n} \\ 
      vel_n = u|_{e} = 0   \qquad &\Rightarrow& \qquad u|_{e} = 0 \\
      \mbox{\textbf{(c)}} \qquad
      vel_t = \frac{u|_{e} + u|_{nne}}{2} = vel_{\bullet}   \qquad &\Rightarrow& \qquad    u|_{nne} = 2 vel_{\bullet} - u|_{e} \\ 
      vel_n = v|_{n} = 0   \qquad &\Rightarrow& \qquad v|_{n} = 0 \\
      \mbox{\textbf{(d)}} \qquad
      vel_t = \frac{u|_{e} + u|_{sse}}{2} = vel_{\bullet}   \qquad &\Rightarrow& \qquad    u|_{sse} = 2 vel_{\bullet} - u_{e} \\ 
      vel_n = v|_{s} = 0   \qquad &\Rightarrow& \qquad v|_{s} = 0
    \end{eqnarray*}
\newpage
\item[Case 3:] The prescribed injection condition (CIPR) has the form $vel_t = 0 \mbox{ and } \linebreak vel_n = vel_I$
    \begin{eqnarray*}
      \mbox{\textbf{(a)}} \qquad
      vel_t = \frac{v|_{nww} + v|_{n}}{2} = 0   \qquad &\Rightarrow& \qquad  v|_{nww} = - v|_{n} \\ 
      vel_n = u|_{w} = vel_I   \qquad &\Rightarrow& \qquad u|_{w} = vel_I \\
    \mbox{\textbf{(b)}} \qquad
      vel_t = \frac{v|_{n} + v|_{nee}}{2} = 0   \qquad &\Rightarrow& \qquad    v|_{nee} = - v|_{n} \\ 
      vel_n = u|_{e} = vel_I   \qquad &\Rightarrow& \qquad u|_{e} = vel_I \\
      \mbox{\textbf{(c)}} \qquad
      vel_t = \frac{u|_{e} + u|_{nne}}{2} = 0   \qquad &\Rightarrow& \qquad    u|_{nne} = - u|_{e} \\ 
      vel_n = v|_{n} = vel_I   \qquad &\Rightarrow& \qquad v|_{n} = vel_I \\
      \mbox{\textbf{(d)}} \qquad
      vel_t = \frac{u|_{e} + u|_{sse}}{2} = 0   \qquad &\Rightarrow& \qquad    u|_{sse} = - u|_{e} \\ 
      vel_n = v|_{s} = vel_I  \qquad &\Rightarrow& \qquad v|_{s} = vel_I
   \end{eqnarray*}

\item[Case 4:] The continuous ejection condition (CECO) is indicated by \linebreak
$\displaystyle \frac{\partial{}}{\partial{n}} vel_t = 0 \mbox{ and } \frac{\partial{}}{\partial{n}} vel_n = 0$
    \begin{eqnarray*}
    \mbox{\textbf{(a)}} \qquad
      \frac{\partial{}}{\partial{n}}(vel_t) = \frac{v|_{nww} - v|_{n}}{\Delta \xi} = 0   \qquad &\Rightarrow& \qquad    v|_{nww} = v|_{n} \\ 
      \frac{\partial{}}{\partial{n}}(vel_n) = \frac{u|_{e} - u|_{w}}{\Delta \eta} = 0   \qquad &\Rightarrow& \qquad    u|_{e} = u|_{w} \\
      \mbox{\textbf{(b)}} \qquad
      \frac{\partial{}}{\partial{n}}(vel_t) = \frac{v|_{nee} - v|_{n}}{\Delta \xi} = 0   \qquad &\Rightarrow& \qquad    v|_{nee} = v|_{n} \\ 
      \frac{\partial{}}{\partial{n}}(vel_n) = \frac{u|_{eee} - u|_{e}}{\Delta \eta} = 0   \qquad &\Rightarrow& \qquad    u|_{eee} = u|_{w} \\
      \mbox{\textbf{(c)}} \qquad
      \frac{\partial{}}{\partial{n}}(vel_t) = \frac{u|_{nne} - u|_{e}}{\Delta \eta} = 0   \qquad &\Rightarrow& \qquad    u|_{nne} = u|_{e} \\ 
      \frac{\partial{}}{\partial{n}}(vel_n) = \frac{v|_{nnn} - v|_{n}}{\Delta \xi} = 0   \qquad &\Rightarrow& \qquad    v|_{nnn} = v|_{n} \\
      \mbox{\textbf{(d)}} \qquad
      \frac{\partial{}}{\partial{n}}(vel_t) = \frac{u|_{sse} - u|_{e}}{\Delta \eta} = 0   \qquad &\Rightarrow& \qquad    u|_{sse} = u|_{e} \\ 
      \frac{\partial{}}{\partial{n}}(vel_n) = \frac{v|_{n} - v|_{s}}{\Delta \xi} = 0   \qquad &\Rightarrow& \qquad    v|_{n} = v|_{s} \\
   \end{eqnarray*}
\end{description}

\section{Numerical results}

In this section we present the numerical results obtained from the proposed method in the previous sections.
The first 
problem relates to the study of the flow between two parallel plates, the second deals with the flow in a square cavity with upper wall 
moving and the third refers to atherosclerosis.
The mesh, for each problem studied, it was generated exactly as described in \cite{cirilo06,saita17}. The edges were built by Spline, and the 
mesh via numerical solution of the grid generation equations.
For mesh generation has been stipulated a maximum of 1000 iterations to solve the linear system that create the grid, considering 
an error of less than $10^{-4}$ in this system. The numerical method of Gauss-Seidel was used to solve the linear system. 

With the first problem we show that our code is able to obtain the numerical solution accurately when compared with analytical
solution.
The second case in our study aims at demonstrate that even for Reynolds number slightly high ($Re = 1000$), where the term convective 
is dominant, it can properly address the problem.
Besides this, with previous cases, show that the proposed numerical model simulates cases whose geometry is cartesian.
Finally, the third case, detailing the simulation when the computational mesh is perfect adjusted the geometry of the problem. 
In this case, even with the dominant convective term ($Re = 900$), the results are in accordance with literature too.

The color maps used in this study range from dark blue (lower velocity) to dark red (higher velocity),
as shown in the Fig. \ref{faixa} below. The results obtained for these problems are presented in the following.
\begin{figure}[!ht]
\begin{center}
\caption{Color maps to the velocity}
\includegraphics[scale=0.25]{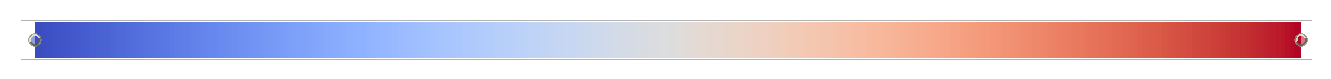}
%\caption{Color maps to the velocity}
\label{faixa}
\center{{\bf{Source:}} The Authors}
\end{center}
\end{figure}

\subsection{Parallel plates problem}

\begin{figure}[!ht]
\begin{center}
\caption{Dimensions of the geometry of the problem laminar flow between two parallel plates \cite{bono11}}
\includegraphics[scale=0.5]{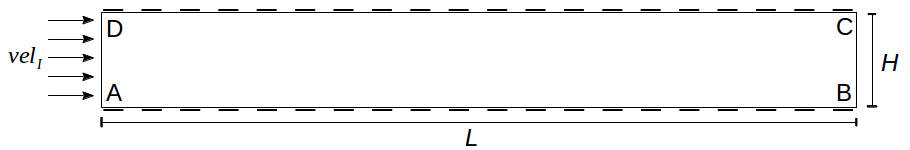}
%\caption{Dimensions of the geometry of the problem laminar flow between two parallel plates \cite{bono11}}
\label{dimensoespp}
\center{{\bf{Source:}} The Authors}
\end{center}
\end{figure}

In this problem we consider a geometry in rectangular shape with height $H$ and length $L = 8H$, according Fig. \ref{dimensoespp}.
The fluid ($Re = 100$) is injected into the geometry by edge $AD$, where the boundary condition is CIPR and  
$vel_I = 1.0 $ $m/s$. The plates are indicated by edges $AB$ and $CD$, where the CNEI conditions is applied.
The fluid output occurs by the edge $BC$, where the boundary condition is CECO.
For the study of the problem we performed five simulations, taking $H = 1$ $m$ and varying the amount of lines considered in the 
directions $\xi$ and $\eta$, as shown in the table \ref{tab:simulations}.

\begin{table}[] 
\begin{center}
\caption{Meshes for simulations to first problem.}
\begin{tabular}{|c|c|c|}
\hline
Mesh & Number of lines in  & Number of lines in   \\
     & the direction $\xi$ & the direction $\eta$  \\ \hline
$P1$ & 9 & 5  \\ \hline
$P2$ & 17 & 9  \\ \hline
$P3$ & 33 & 17   \\ \hline
$P4$ & 65 & 33  \\ \hline
$P5$ & 129 & 65  \\ \hline
\end{tabular}
%\caption{Meshes for simulations to first problem.}
\label{tab:simulations}
\center{{\bf{Source:}} The Authors}
\end{center}
\end{table}
The initial condition was taken from the state of quiescence. The simulation was performed until the steady state was reached, and 
this occurred at $ \tau = 30$, approximately.
In this study, for convergence of the pressure and velocity equations we consider the 
$\Delta \tau =10^{-2}$ to meshes $P1$ to $P3$. 
But to $P4$ and $P5$ the $\Delta \tau$ taken was equal to $5.10^{-3}$.
This difference values in $ \Delta \tau $ occurred because of the mesh refinement.

\begin{figure}[!ht]
\begin{center}
\caption{ Velocity profile in the output section to flow between two parallel plates for different mesh refinements}
\subfigure[In this paper]{\includegraphics[scale=0.23]{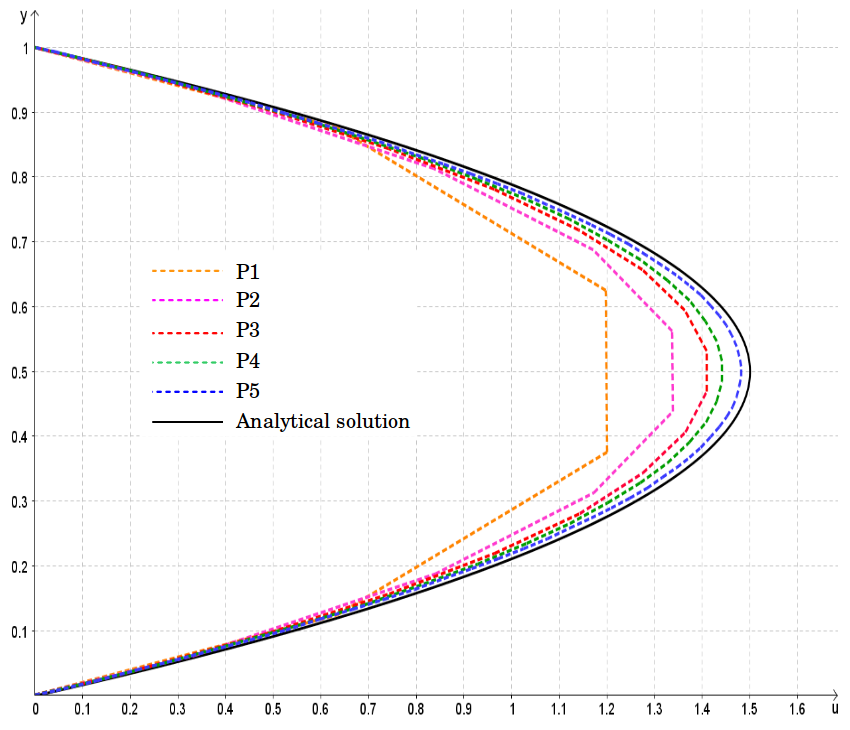}
  \label{grafpp}}\qquad \qquad
\subfigure[Presented by \cite{bono11}]{\includegraphics[scale=0.23]{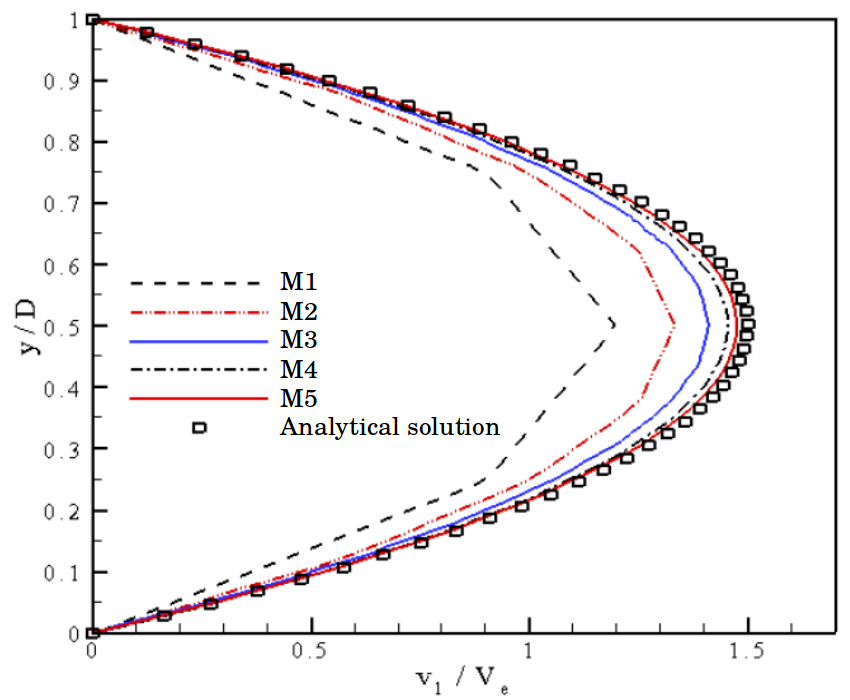}
   \label{grafppref}} 
%\caption{ Velocity profile in the output section to flow between two parallel plates for different mesh refinements.}
\label{grafpparalelas}  
\center{{\bf{Source:}} The Authors}
\end{center}
\end{figure}

In Figs. \ref{grafpp} and \ref{grafppref} the velocity profiles in the output section, in the edge $BC$, are presented \footnote{In this 
case we are only analyzing component $u$ of the velocity because this is perpendicular to the edge where there is fluid output.} for the 
meshes considered in this work and in \cite{bono11}, respectively. 
Note that in \cite{bono11} $D = 1$, then $y/D = y$.
Likewise, as $v_e = 1$ follows $v_l/v_e = v_l = u$.
Therefore, we can compare the graphs shown in the figures because both describe the same relationship.

Looking at the Fig. \ref{grafpp}, can be observed that further refinement of the mesh gives better results.
In the mesh $P5$, for example, the maximum value found for the speed was equal to $1.4814$ $m/s$.
When comparing the profiles obtained by the proposed method in this work and literature \cite {bono11}, note that both are close 
to the analytical solution. However, in the literature the maximum value was 1.4756 m/s with the finite elements method.
It is known in the scientific community which computational cost of finite differences (our proposal) used is
lower than finite elements.
   \begin{table}[] 
     \centering
		  \caption {Convergence speed of the numerical code via free stream ${V_l}_{num}$}
	    \begin{tabular}{|c|c|c|c|c|c|}
	      \hline
	      Meshes		& ${V_l}_{num}$ 		& $h_x$			& $h_y$			& $h_x \times h_y$		& $Error=|1.5 - {V_l}_{num}|$ \\
	      \hline 
	      $P_1$		& \small{1.1999} 	& \small{1.0}		& \small{0.25}		& \small{0.25}			& \small{0.3001}	\\
	      \hline 
	      $P_2$		& \small{1.3380} 	& \small{0.5}		& \small{0.125}		& \small{0.0625}			& \small{0.1620}	\\
	      \hline 
	      $P_3$		& \small{1.4158} 	& \small{0.25}		& \small{0.0625}		& \small{0.015625}		& \small{0.0842}	\\
	      \hline 
	      $P_4$		& \small{1.4498} 	& \small{0.125}		& \small{0.03125}	& \small{0.00390625}		& \small{0.0502}	\\
	      \hline 
	      $P_5$		& \small{1.5068} 	& \small{0.0625}		& \small{0.015625}	& \small{0.0009765625}		& \small{0.0068}	\\
	      \hline
	    \end{tabular}
     %\caption {Convergence speed of the numerical code via free stream ${V_l}_{num}$}
     \label{tab:erro}
		\center{{\bf{Source:}} The Authors}
 \end{table}
In Table \ref{tab:erro} is shown that the mesh refinement implies lower error between the theoretical velocity and the 
calculated velocity in our code. So there is clearly a process of convergence to the theoretical velocity, i.e., the mesh doesn't 
interfere in the solution obtained. Finally, the velocity field obtained, when the steady state has been reached, can be 
seen in the Fig. \ref{grafvelocidade} below.

\begin{figure}[!ht]
\begin{center}
\caption{Velocity field obtained in this paper to first problem}
\includegraphics[scale=0.5]{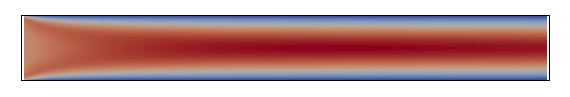}
%\caption{Velocity field obtained in this paper to first problem}
\label{grafvelocidade}
\center{{\bf{Source:}} The Authors}
\end{center}
\end{figure}

\subsection{Lid-driven cavity problem}

\begin{figure}[!ht]
\begin{center}
\caption{Geometry and dimensions considered in the problem of square cavity} 
\includegraphics[scale=0.45]{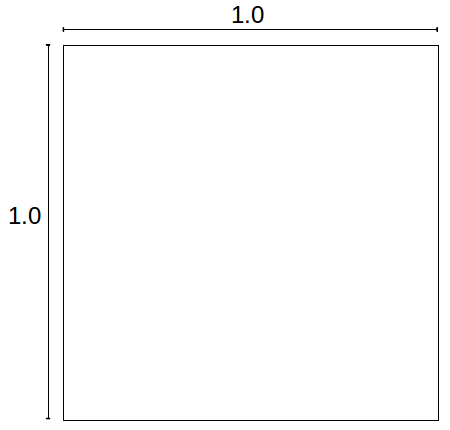}
%\caption{Geometry and dimensions considered in the problem of square cavity} 
\label{cavidadedim}
\center{{\bf{Source:}} The Authors}
\end{center}
\end{figure}

In this problem we consider a square cavity with edge measuring $1$ $m$, whose geometry is shown in Fig. \ref{cavidadedim}.
It is completely fluid filled, the upper wall is moving and we apply the boundary condition CLES with slip velocity 
equal to $1 m/s$. The others walls are subject to boundary conditions CNEI.

The aim of this problem is the numerical verification of algorithm.
The convective terms are more evident and the numerical model is more required numerically.
The results obtained for this problem are compared with those presented by \cite{bono11,ghia82,griebel98,gupta05,hou95,marchi09}.

This study was carried out from a mesh with 129 lines $\xi$ and $\eta$.
The initial condition was taken from the state of quiescence. The simulation was performed until the steady state was reached, and 
this occurred at $ \tau = 50$, approximately.
For convergence of the pressure and velocity equations, we consider the $\Delta \tau =10^{-3}$ because of the mesh refinement.

\begin{figure}[!ht]
\begin{center}
\caption{ Velocity fields obtained in solving the problem of the cavity in steady state}    
\subfigure[$Re = 100$ ]{\includegraphics[scale=0.297]{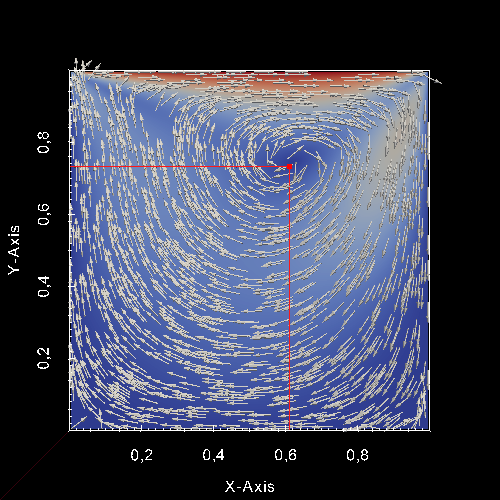}
  \label{6cav100vec}}\qquad \quad \quad
\subfigure[$Re = 400$ ]{\includegraphics[scale=0.29]{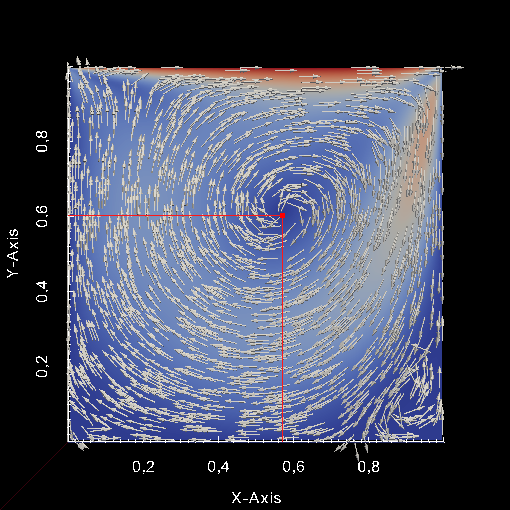}
  \label{6cav400vec}} 
%\caption{ Velocity fields obtained in solving the problem of the cavity in steady state}    
\label{cavperfis}
\center{{\bf{Source:}} The Authors}
\end{center}
\end{figure}

The variations of velocity fields obtained, after steady state has been reached, can be seen in the 
Figs. \ref{6cav100vec} and \ref{6cav400vec}.
In each of the figures we can see the 
formation of a primary vortex near the center of the domain, whose location varies with Reynolds number considered. 
Moreover, we can observe the formation of two other smaller vortices in the lower regions near to 
left and right boundaries of the cavity.
In the first case, where $Re = 100$, there is only an indication of the formation of secondary vortices.
For $Re = 400$ we observe a vortex formed in the lower right corner, and only one beginning of the vortex in the left corner. 
It is due to increased Reynolds number.

From the localization of the vortices, the interest was to compare the main vortex coordinates obtained with 
our methodology in relation to other studies.
The table \ref{tab:Comparações} present these coordinates for cases where $Re = 100$ and $Re = 400$.
Note that our results are in agreement with the literature.

 \begin{table}[] 
     \centering
		     \caption {Localization of the main vortex center}
   \begin{tabular}{|c|c|c|}
     \hline
	Reference & $Re = 100$ & $Re = 400$  \\ \hline
	This work & $(0.6109, 0.7335)$ & $(0.5699, 0.6033)$  \\ \hline
	Bono et al. \cite{bono11} & $(0.6157, 0.7373)$ & $(0.5613, 0.6123)$  \\ \hline
	Ghia et al. \cite{ghia82} & $(0.6172, 0.7344)$ & $(0.5547, 0.6055)$   \\ \hline
	Gupta e Kalita \cite{gupta05} & $(0.6125, 0.7375)$ & $(0.5500, 0.6125)$  \\ \hline
	Hou et al. \cite{hou95} & $(0.6196, 0.7373)$ & $(0.5608, 0.6078)$  \\ \hline
	Marchi et al. \cite{marchi09} & $(0.6162, 0.7373)$ & $(0.5537, 0.6054)$ \\ \hline
   \end{tabular}
     %\caption {Localization of the main vortex center}
     \label{tab:Comparações}
 \center{{\bf{Source:}} The Authors}
 \end{table}

One third case ($Re = 1000$), the secondary 
vortices are very defined in the lower right and left corners. There is evidence of the formation of another vortex in the
upper left corner. Analogously, it is due to increased Reynolds number too. 

\begin{figure}[!ht]
\begin{center}
\caption{ Velocity fields to second problem indicating the direction of flow to $Re = 1000$}    
\subfigure[Profile obtained in this work. {\bf{Source:}} The Authors]{\includegraphics[scale=0.25]{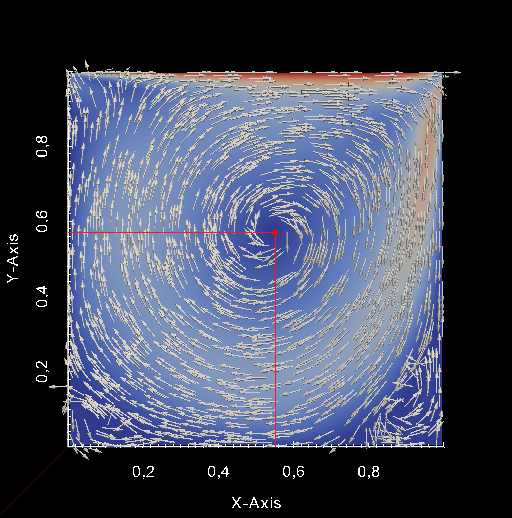}
  \label{6cav1000vecx}}\qquad \quad
\subfigure[Profile presented by \cite{griebel98}]{\includegraphics[scale=0.8]{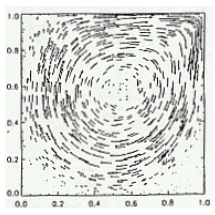}
  \label{6cav1000vecref}}
%\caption{ Velocity fields to second problem indicating the direction of flow to $Re = 1000$}    
\label{cavperfil1000}
\end{center}
\end{figure}

Graphically, the results displayed in Fig. \ref{6cav1000vecx}, compared with that obtained in Fig. \ref{6cav1000vecref}, shows that our simulation is in accordance with literature \cite{griebel98}. 
The coordinates of the primary vortex obtained with our numerical code were (0.5528, 0.5698).

\subsection{Atherosclerosis problem}

Atherosclerosis is a disease associated with accumulation of lipids, complex carbohydrates, blood components, cells and other elements 
in large and medium-sized arteries, it is the main cause of heart disease \cite{souza05}. In general, the development of this 
problem starts from accumulation of cholesterol LDL in the artery walls type, which may be higher or lower depending on the availability 
of this substance in the blood \cite{lusis00}. The accumulation of these compounds in the walls of an artery causes the same 
hardening through the formation of atherosclerotic plaques, which can lead to stenosis, or narrowing of the blood vessel, reducing 
blood flow in the artery \cite{fukujima99}.

\begin{figure}[!ht]
\begin{center}
\caption{Geometry and dimensions considered on the atherosclerosis problem}
\includegraphics[scale=0.4]{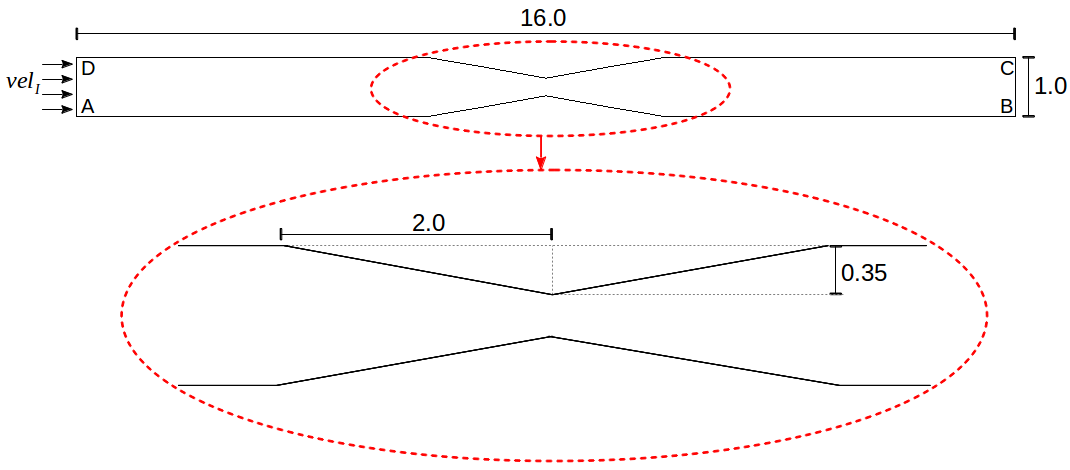}
%\caption{Geometry and dimensions considered on the atherosclerosis problem}
\label{artdimens}
\center{{\bf{Source:}} The Authors}
\end{center}
\end{figure}

Motivated by this problem, our objective was to reproduce a single case of blood flow in the region of a large caliber artery 
containing a stenosis.
The constriction in the upper and lower border, both with same size, is located in the same point as 
presented in Fig. \ref{artdimens}.

\begin{figure}[!ht]
\begin{center}
\caption{ Mesh considered in solving the problem related to atherosclerosis}
\includegraphics[scale=0.42]{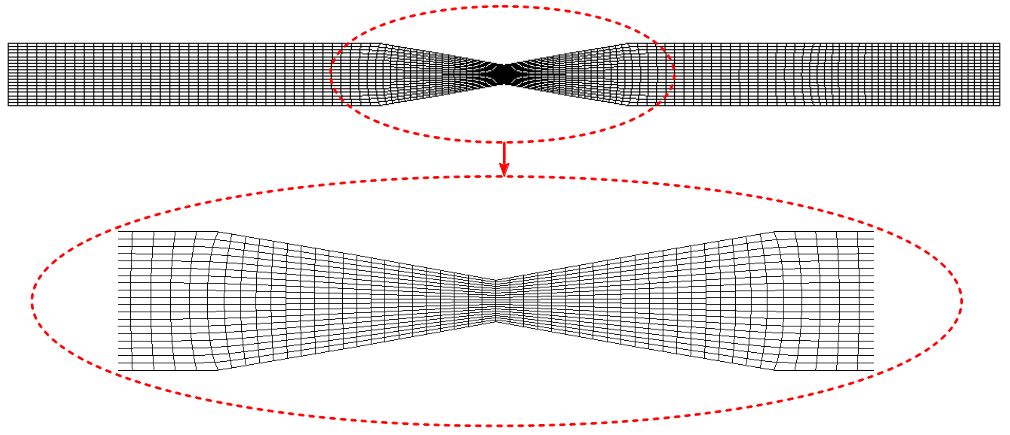}
%\caption{ Mesh considered in solving the problem related to atherosclerosis}
\label{artmalha}
\center{{\bf{Source:}} The Authors}
\end{center}
\end{figure}

We consider a fluid with $Re = 900$. This consideration is acceptable to the scientific community, and we can approach the 
blood as a Newtonian incompressible viscous fluid \cite{layek07}.
The fluid is injected into the geometry by edge $AD$, subject to the boundary condition CIPR with $vel_I = 0.1467$.
The output occurs through edge $BC$, where applied condition is CECO, and in the other walls we consider the boundary condition CNEI.
As in previous cases, the initial condition was taken from the state of quiescence.
For the study we consider the mesh shown in Fig. \ref{artmalha}, constructed from the dimensions set 
out above, containing $\xi \times \eta = 129 \times 20$ lines, and $\Delta \tau = 5.10^{-3}$.

\begin{figure}[!ht]
\begin{center}
\caption{ Velocity field obtained from the simulation of the problem related to atherosclerosis}
\includegraphics[scale=0.4]{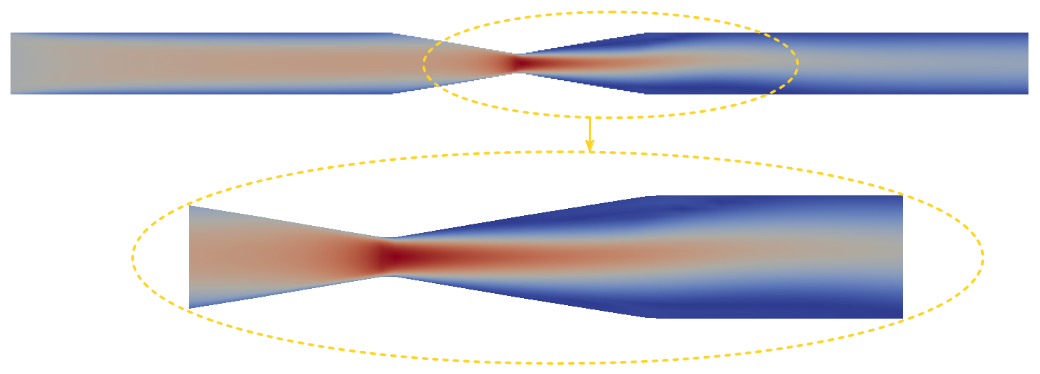}
%\caption{ Velocity field obtained from the simulation of the problem related to atherosclerosis}
\label{artvelocid}
\center{{\bf{Source:}} The Authors}
\end{center}
\end{figure}

The velocity field obtained after steady state is reached is shown in Fig. \ref{artvelocid}.
The fluid enters in the geometry, suffers the performance of CNEI boundary condition, and the fully developed profile is established.
When the fluid begins its entry path by stenosis undergoes a continuous drop in pressure since the Reynolds is moderated, 
thereby increasing the velocity. The maximum velocity value occurs in the narrowing region.
After passing through narrowing the fluid velocity decreases and flows throughout the rest of the geometry, 
but in a lower level than that of the input.
This pattern in the flow is also observed by \cite{fukujima99}, which leads to implications for human health.

\begin{figure}[!ht]
\begin{center}
\caption{Velocity field for the components $u$ (left) and $v$ (right) in the forming region of the vortexes}
\includegraphics[scale=0.19]{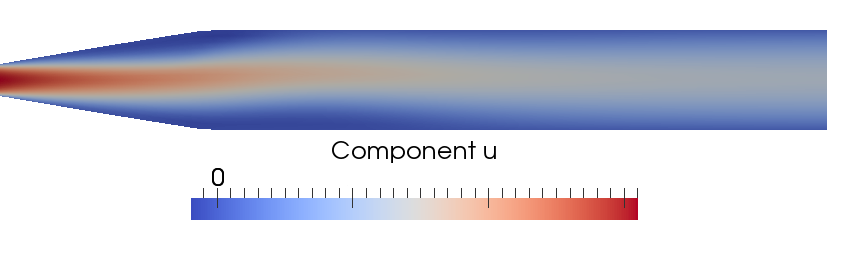}
\quad
\includegraphics[scale=0.19]{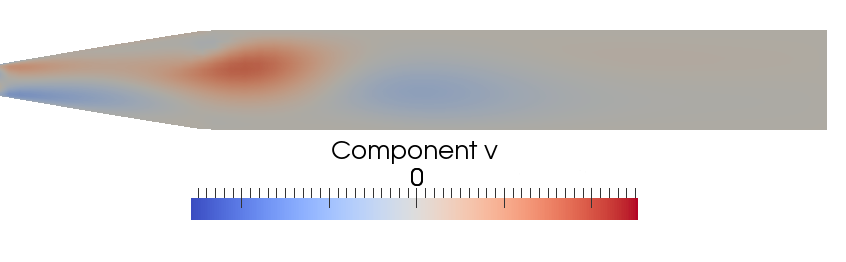}
%\caption{Velocity field for the components $u$ (left) and $v$ (right) in the forming region of the vortexes}
\label{artvelocinfsup}
\center{{\bf{Source:}} The Authors}
\end{center}
\end{figure}

In Fig. \ref{artvelocinfsup} are plotted the values of $u$ and $v$ components, respectively.
From these graphs we can see that after the stenosis there is occurrence of vortex formation.
Note that although the geometry has symmetry, the velocity components are asymmetric in the domain.
The vortex near the top wall, in size and position with respect to the abscissa, is different from the vortex that is
side of the bottom wall.
This asymmetry in the velocity field is due to the moderate value of the Reynolds number used in the simulation.
This asymmetric behavior was also obtained by \cite{layek07}.
For Reynolds values less than $900$ the flow tends to be symmetrical, because the vortices tend to disappear.
But for values greater than $900$ the convective terms are strongly dominant, the flow tends to be asymmetric and
vortices increase in size, intensity and position with respect to abscissa.
The simulation here displayed when $Re = 900$ is a turbulent presage.

The study of the disordered formation of these vortices beyond the stenosis is critical \cite{layek07}, once they can hinder the continuity of blood flow in arteries and thus, can aggravate health problems that may lead the patient to death.

\section{Conclusion}

The objective of this study was to present an effective numerical model to solve the equations of fluid dynamics
in the case: laminar, newtonian, incompressible, isothermal, two-dimensional in curvilinear coordinates system.
In addition, show that even applying a first order scheme (FOU) in the convective terms, it was possible to obtain 
satisfactory results with the numerical model.

In the cases of flows between parallel plates and in the lid-driven cavity, the numeric code was able to simulate the velocity field accurately,
that is, the numerical model converged to values of interest does not depend on the mesh (see table \ref{tab:erro}) and also for 
various values of Reynolds (see Figs. \ref{grafpparalelas}, \ref{cavperfis}, \ref{cavperfil1000}).
Furthermore, we show that in the case of geometries whose cartesian coordinate system is completely inappropriate, the  
curvilinear coordinates system appears as a prominent alternative. This was the case of atherosclerosis problem.
In this problem we show that the mesh is perfectly adequate to the domain, and that the proposed numerical model allowed to obtain 
numerical solutions in agreement with the literature.

For future work we intend to extend the methodology adopted in this study. 
First we will implement a higher-order convective scheme which solves problems in a more comprehensive 
range of Reynolds number, and then extend the code for three-dimensional cases.


\begin{thebibliography}{8}
\providecommand{\url}[1]{{#1}}
\providecommand{\urlprefix}{URL }
\expandafter\ifx\csname urlstyle\endcsname\relax
  \providecommand{\doi}[1]{DOI~\discretionary{}{}{}#1}\else
  \providecommand{\doi}{DOI~\discretionary{}{}{}\begingroup
  \urlstyle{rm}\Url}\fi

\bibitem{Amsden_1970}
AMSDEN, A.A.; HARLOW, F.H. A simplified MAC technique for incompressible fluid flow calculations. {\it Journal of Computational Physics}, v. 6, p. 322-325, 1970.

\bibitem{bono11}
BONO, G.; LYRA, P.R.M.; BONO, G.F.F. Solu\c{c}\~ao num\'erica de escoamentos
incompress\'{\i}veis com simula\c{c}\~ao de grandes escalas. {\it Mec\'anica Computacional}, v. 30, p. 1423-1440, 2011.

\bibitem{cirilo06}
CIRILO, E.R.; DE BORTOLI, A.L. Cubic splines for trachea and bronchial tubes grid generation. {\it Semina: Exact and Technological Sciences}, v. 27, p.147-155, 2006.

\bibitem{courant52}
COURANT, R.; ISAACSON, E.; REES, M. On the solution of nonlinear hyperbolic differential equations by finite difference. 
{\em Communications on Pure and Applied Mathematics}, v. 5, p. 243-255, 1952.

\bibitem{ferreira12}
FERREIRA, V.G.; LIMA, G.A.B.; CORREA, L.; CANDEZANO, M.A.C.; CIRILO, E.R.; ROMEIRO, N.M.L.; NATTI, P.L. Avalia\c{c}\~ao computacional de esquemas convectivos em problemas de din\^amica dos fluidos. {\it Semina: Exact and Technological Sciences}, v. 32, p. 107-116, 2012.

\bibitem{fortuna12}
FORTUNA, A.O. Computational techniques for fluid dynamics: Basic concepts and applications. S\~ao Paulo: EDUSP, 2012.

\bibitem{fukujima99}
FUKUJIMA, M.M.; GABBAI, A.A. Condutas na estenose da car\'otida. {\it Revista Neuroci\^encias}, v. 7, p. 39-44, 1999.

\bibitem{ghia82}
GHIA, U.; GHIA, K.N.; SHIN, C.T. High-Re solutions for incompressible flow
using the Navier-Stokes equations and multigrid method. {\it Journal of Computational Physics}, v. 48, p. 387-411, 1982.

\bibitem{griebel98}
GRIEBEL, M.; DORNSEIFER, T.; NEUNHOEFFER, T. Numerical simulation in fluid
dynamics: A pratical introduction. Philadelphia: Society for Industrial and Applied Mathematics Press, 1998.

\bibitem{gupta05}
GUPTA, M.M.; KALIT, J.C. A new paradigm for solving Navier-Stokes equations:
streamfunction-velocity formulation. {\it Journal of Computational Physics}, v. 207, p. 52-68, 2005.

\bibitem{hou95}
HOU, S.; ZOU, Q.; CHEN, S.; DOOLEN, G.; COGLEY, A. Simulation of cavity flows
by the lattice Boltzmann method. {\it Journal of Computational Physics}, {\bf 118} (1995) 329--347.

\bibitem{Huang_2013}
HUANG, Y.L.; LIU, J.G.; WANG, W.C. A generalized MAC scheme on curvilinear domains. {\it SIAM Journal on Scientific Computing}, v. 35, p. B953-B986, 2013.

\bibitem{johnson_1987}
JOHNSON, C. Numerical solution of partial differential equations by the 
finite element method. New York: Cambrigde University Press, 1987.

\bibitem{layek07}
LAYEK, G.C.; MIDYA, C. Effect of constriction height on flow separation in a two-dimensional channel. {\it Communications in Nonlinear Science and Numerical Simulation}, v. 12, p. 745-759, 2007.

\bibitem{leonard_1979}
LEONARD, B.P. A stable and accurate convective modelling procedure based on quadratic upstream interpolation. {\it Computer Methods in Applied Mechanics and Engineering}, v. 19, p. 59-98, 1979.

\bibitem{leveque_2002}
LEVEQUE, R.J. Finite volume methods for hyperbolic problems. New York: Cambrigde University Press, 2002.

\bibitem{lusis00}
LUSIS, A.J. Atherosclerosis. {\it Nature}, v. 407, p. 233-241, 2000.

\bibitem{maliska_1983}
MALISKA, C.R.; RAITHBY, G.D. Calculating three-dimensional fluid flows using nonorthogonal grids. In: 3rd International Conference in Numerical Methods in Laminar and Turbulent Flow, San Francisco: Pineridge Press, 1983.
\newblock
\urlprefix\url{http://www.sinmec.ufsc.br/site/arquivos/f-fwhkiqlbnq_1983_3dimensional.pdf}

\bibitem{maliska13}
MALISKA, C.R. Heat transfer and computational fluid mechanics. Rio de Janeiro: LTC, 2013.

\bibitem{marchi09}
MARCHI, C.H.; SUERO, R.; ARAKI, L.K. The lid-driven square cavity flow:
Numerical solution with a 1024x1024 grid. {\it Journal of the Brazilian Society of Mechanical Sciences and Engineering}, v. 31, p. 186-198, 2009.

\bibitem{McKee_2008}
MCKEE, S.; TOM\'E, M.F.; FERREIRA, V.G.; CUMINATO, J.A.; CASTELO, A.; SOUSA, F.S.; MANGIAVACCHI, N. The MAC methods. {\it Computers and Fluids}, v. 37, p. 907-930, 2008.

\bibitem{Monaghan_2005}
MONAGHAN, J.J. Smoothed particle hydrodynamics. {\it Reports on Progress in Physics}, v. 68, p. 1703-1759, 2005.

\bibitem{Pardo_2012}
PARDO, S.R.; NATTI, P.L.; ROMEIRO, N.M.L.; CIRILO, E.R. A transport modeling of the carbon-nitrogen cycle at Igap\'o I Lake-Londrina,
Paran\'a State, Brazil. {\it Acta Scientiarum. Technology}, v. 34, p. 217-226, 2012.

\bibitem{Patil_2009}
PATIL, P.P.; TIWARI, S. Computation of flow past complex geometries using MAC algorithm on body-fitted coordinates. {\it Engineering Applications of Computational Fluid Mechanics}, v. 3, p. 15-27, 2009.

\bibitem{Price1966}
PRICE, H.S.; VARGA, R.S.; WARREN, J.E. Applications of oscillation matrices to diffusion-correction equations. {\it Journal of Mathematical Physics}, v. 45, p. 301-311, 1966.

\bibitem{Romeiro_2011}
ROMEIRO, N.L.M.; CASTRO, R.G.S.; CIRILO, E.R.; NATTI, P.L. Local
calibration of coliforms parameters of water quality problem at Igap\'o I
Lake - Londrina, Paran\'a, Brazil. {\it Ecological Modelling}, v. 222, p. 1888-1896, 2011.

\bibitem{Romeiro17}
ROMEIRO, N.L.M.; MANGILI, F.B.; COSTANZI, R.N.; CIRILO, E.R.; NATTI, P.L. Numerical simulation of BOD5 dynamics in Igap\'o I lake, Londrina, Paran\'a, Brazil: Experimental measurement and mathematical modeling. {\it Semina: Exact and Technological Sciences}, v. 38, p. 50-58, 2017.

\bibitem{saita17}
SAITA, T.M.; NATTI, P.L.; CIRILO, E.R.; ROMEIRO, N.L.M.; CANDEZANO, M.A.C.; ACUNA, R.A.B.; MORENO, L.C.G. Simula\c{c}\~ao num\'erica da din\^amica de coliformes fecais no lago Luruaco, Col\^ombia. {\it Trends in Applied and Computational Mathematics}, v. 18, p. 435-447, 2017.	

\bibitem{santos_2012}
SANTOS, F.L.P.; FERREIRA, V.G.; TOM\'E, M.F.; CASTELO, A.; MANGIAVACCHI, N.; MCKEE, S. A marker-and-cell approach to free surface 2-D multiphase flows. {\it International Journal for Numerical Methods in Fluids}, v. 70, p. 1543-1557, 2012.

\bibitem{Shah_2012}
SHAH, A.; YUAN, L.; ISLAM, S. Numerical solution of unsteady Navier-Stokes equations on curvilinear meshes. {\it Computers and Mathematics with Applications}, v. 63, p. 1548-1556, 2012.

\bibitem{souza05}
DE SOUZA, L.V.; DE CASTRO, C.C.; CERRI, G.G. Evaluation of carotid  
atherosclerosis by ultrasound and magnetic resonance imaging. {\it Radiologia Brasileira}, v. 38, p. 81-94, 2005.

\bibitem{thompson76}
THOMPSON, J.F.; THAMES, F.C.; MASTIN, W.C. Boundary fitted curvilinear
coordinate system for solution of partial differential equations on fields
containing any number of arbitrary two-dimensional bodies. {\it NASA Langley Research Centre}, v. CR-2729, 1976.
\urlprefix\url{http://ntrs.nasa.gov/archive/nasa/casi.ntrs.nasa.gov/19770021145.pdf}

\bibitem{thompson85}
THOMPSON, J.F.; WARSI, Z.U.A.; MASTIN, C.W. Numerical grid generation: Foundations and applications. New York: Elsevier North-Holland, 1985.

\bibitem{tome_2014}
TOM\'E, M.F.; CASTELO, A.; N\'OBREGA, J.M.; PAULO, G.S.; PEREIRA, F.T. Numerical and experimental investigations of three-dimensional container filling with Newtonian viscous fluids. {\it Computers and Fluids}, v. 90, p. 172-185, 2014.


\end{thebibliography}
\end{document}